\documentclass[a4paper,12pt]{article} 
\usepackage{graphicx,latexsym,amsmath,amsthm,amssymb,color,indentfirst,xfrac}
 \usepackage[usenames,dvipsnames]{pstricks}
\usepackage{mathrsfs} %ตัวเขียน
 \usepackage{epsfig}
\usepackage{pst-grad} % For gradients
 \usepackage{pst-plot} % For axes

\newtheorem{satz}{thm1}

\newtheorem{sat}{thm2}

\newtheorem{thm}[sat]{Theorem}
\newtheorem{lem}[satz]{Lemma}

\newtheorem{prop}[satz]{Proposition}
\newtheorem{ex}[satz]{Example}
\theoremstyle{definition}

\newcommand{\be}{\begin{equation}}           % equation
\newcommand{\ee}{\end{equation}}

\newcommand{\ba}{\begin{align}}                % align
\newcommand{\ea}{\end{align}}

\newcommand{\bal}{\begin{align*}}              % align*
\newcommand{\eal}{\end{align*}}

\newcommand{\bxx}{\begin{ex}}

\newcommand{\exx}{\end{ex}}

\newcommand{\txsit}

\newenvironment{pr}
{\begin{trivlist}
\item[\hskip\labelsep{\bf Proof.}]}                     
{$\hfill\Box$\end{trivlist}}

\title{Sufficient conditions on cycles that make planar graphs
4-choosable }
\author {Pongpat Sittitrai\\ 
{\small\em Department of Mathematics, Faculty of Science, Khon Kaen University, 40002, Thailand }\\  
{\small\em E-mail address: pongpat\_s@kkumail.com} 
\and Kittikorn Nakprasit \footnote{Corresponding Author} \\ 
{\small\em Department of Mathematics, Faculty of Science, Khon Kaen University, 40002, Thailand }\\
{\small\em E-mail address: kitnak@hotmail.com}}
\date{}

\begin{document}

\maketitle

\begin{center}{\bf Abstract}\end{center}
\indent\indent
Xu and Wu proved that if every $5$-cycle of a planar graph $G$ 
is not simultaneously adjacent to $3$-cycles and $4$-cycles, 
then $G$ is $4$-choosable. In this paper, we improve this result as follows. 
Let $\{i, j, k, l\} = \{3,4,5,6\}.$ For any chosen $i,$ if every $i$-cycle 
of a planar graph $G$ is not simultaneously adjacent to $j$-cycles, 
$k$-cycles, and $l$-cycles, then $G$ is $4$-choosable.

\section{Introduction}
\indent Every graph in this paper is finite, simple, and undirected graph. 
The concept of choosability was introduced by Vizing in 1976 \cite{Vizing} 
and Erd\H os, Rubin, and Taylor in 1979 \cite{Erdos}, independently. 
A $k$-\textit{list assignment} $L$ of a graph $G$ assigns a list $L(v)$ 
(a set of colors) and $|L(v)|= k$ to each vertex $v.$ 
A graph $G$ is $L$-colorable if there is a proper coloring $f$ where $f(v)\in L(v).$   
If $G$ is $L$-colorable for any $k$-assignment $L,$ then we say  $G$ 
is $k$-\textit{choosable}. 

It is known that every planar graphs is $4$-colorable \cite{app1, app2}.  
Thomassen \cite{Tho} proved that every planar graph is $5$-choosable. 
In contrast, Voight \cite{Vo1} presented an example of non $4$-choosable planar graph. 
Additionally, Gutner \cite{Gut} showed that determining 
whether a given planar graph $4$-choosable is NP-hard. 
Since every planar graph without $3$-cycle always has a vertex of 
degree at most $3,$ it is $4$-choosable. 
More  conditions for  a planar graph to be $4$-choosable are investigated. 
It is shown that a planar graph is $4$-choosable if it has no  
$4$-cycles \cite{Lam2}, $5$-cycles \cite{Wang1}, $6$-cycles \cite{Fi}, 
$7$-cycles \cite{Far}, intersecting $3$-cycles \cite{Wang2}, 
intersecting $5$-cycles \cite{Hu},  or 
$3$-cycles adjacent to $4$-cycles \cite{Bo, Cheng}. 
Xu and Wu  \cite{Xu} proved that if every $5$-cycle of a planar graph $G$ 
is not simultaneously adjacent to $3$-cycles and $4$-cycles, 
then $G$ is $4$-choosable. In this paper, we improve this result as follows.

\begin{thm}\label{main2}
Let $\{i, j, k, l\} = \{3,4,5,6\}.$ For any chosen $i,$ if every $i$-cycle  
of a planar graph $G$ is not simultaneously adjacent to $j$-cycles, 
$k$-cycles, and $l$-cycles, then $G$ is $4$-choosable.   
\end{thm}

\section{Structure}
\indent First, we introduce some notations and definitions. 
A $k$-vertex (face) is a vertex (face) of degree $k,$ 
a $k^+$-vertex (face) is a vertex (face) of degree at least $k,$ 
and a $k^-$-vertex (face) is a vertex (face) of degree at most $k.$ 
A \textit{$(d_1,d_2,\dots,d_k)$-face} $f$ is a face of degree $k$ 
where all vertices on $f$ have degree $d_1,d_2,\dots,d_k$. 
A \textit{$(d_1,d_2,\dots,d_k)$-vertex} $v$ is a vertex of degree $k$ 
where all faces incident to $v$ have degree $d_1,d_2,\dots,d_k$. 
A \textit{wheel graph}  $W_n$ is an $n$-vertex graph formed 
by connecting a single vertex (\textit{hub}) to all vertices 
(\textit{external vertices}) of an $(n-1)$-cycle.

Some basic properties are collected in the following proposition.  
\begin{prop}\label{main1} 
Let $\{i, j, k, l\} = \{3,4,5,6\}$ and $G$ be a planar graph 
such that every $i$-cycle is not simultaneously adjacent to $j$-cycles, 
$k$-cycles, and $l$-cycles. Embed $G$ into the plane, then\\  
(1)  $G$ does not contains a $4$-cycle with one chord that shares exactly one edge 
in a cycle with a $4$-cycle or a $5$-cycle.\\
(2)  $G$ does not contains a $4$-face that  shares  exactly one edge simultaneously 
 with a $3$-cycle and a $4^-$-cycle.\\ 
(3)  $G$ does not contains $5$-cycle with one chord that shares exactly one edge 
with two $5$-cycles.\\ 
(4) $G$ does not contains $W_5$ that shares exactly one edge with a $6^-$-cycle. 
\end{prop}

A $(3,3,5,5^+)$-vertex is called a \textit{flaw $4$-vertex}.  
a $5$-face is called  a \textit{poor $5$-face} if it is adjacent to 
at least four $3$-faces, 
and incident with either five $4$-vertices or four $4$-vertices and one $5$-vertex. 
Let $P$ be a face and $T$ a $3$-face in $G$. If $P$ and $T$ share exactly one edge. 
The vertex $v$ on $T$ but not on $P$ is called a \textit{source}.\\

In this part, we consider a minimal non $4$-choosable planar graph $G$ embeded into the plane.

%\begin{figure}[ht]\label{fig2}
%\centering
%\includegraphics[scale=.5]{donthave2.eps}
%\caption{These graphs follow the conditions of Proposition \ref{main1}}
%\end{figure}

\begin{lem}\cite{Lam2} Every vertex has degree at least $4.$ 
\end{lem}

\begin{lem}\cite{Bo}\label{poor} Every source of poor $5$-face is a $5^+$-vertex.
\end{lem}
\begin{lem}\cite{Bo}\label{lema} Let $f$ and $g$ be two faces in $G$. 
If $f$ shares exactly one edge $x_i{x_j}$ with $g$ where $d(x_i)\leq 5$, 
then at least one vertex in $V(f)\cup V(g)-\{x_i\}$ is a $5^+$-vertex. 
\end{lem}
\begin{thm}\label{thma} Let $v$ be a $5$-vertex in $G$. Then at least three incident faces of $v$ are incident to at least two $5^+$-vertices.
\end{thm}
\begin{pr} Suppose that $v$ has  at most two incident faces that are incident 
to at least two $5^+$-vertices. Then there are two incident faces 
of $v,$ $f_1$ and $f_2$ where $f_1$ and $f_2$ are adjacent and are   
$(4,4,\dots,4,5)$-faces, a contradiction to Lemma \ref{lema}.
\end{pr}
\begin{thm}\label{thmb} Let $v$ be a $6$-vertex in $G$. Then at least two  incident faces of $v$  are incident at least two $5^+$-vertices.
\end{thm}

\begin{pr}  Let $v$ be incident to six faces $f_1$, $f_2,\dots, f_6$. 
Suppose that at most one $f_i$ is a $(4^+,4^+,\dots,5^+,6)$-face. 
By minimality of $G$ and for any subgraph $Z$, the graph $G -V(Z)$ has an $L'$-coloring. If we show that $Z$ has an $L''$-coloring for any  $|L''(x_i)|=4-|N(x_i)-Z|)$, then we obtain $G$ is $4$-choosable, a contradiction. 
Let each $f_i$ be a $(4,4,\dots,4,6)$-face and $Z$ be a graph $f_1\cup f_2\cup\dots\cup f_6$. Then $L''(v)=4$, $L''(x)=3$ if $x$ 
is adjacent to $v$, and $L''(x)=2$ for otherwise. We choose one adjacent vertex of $v$ given $w$.  There is a color $a$ in $L''(v)-L''(w)$ and we color $v$ with a color $a$. Moreover, $|L''(w)-\{a\}|=3$ and $|L''(x)-\{a\}|\geq2$ for each $x$. Thus $Z$ has an $L''$-coloring since each cycle is $2$-choosable if there are two lists that are not equal, a contradiction.
Let $v$ be incident exactly one $(4^+,4^+,\dots,5^+,6)$-face given $f_6$,  $Z$ be a graph $f_1\cup f_2\cup\dots\cup f_5$, and  $x_1$ and $x_2$ be 
two vertices that adjacent to $v$ and  incident to $f_6.$ 
Then $|L''(x)|=2$ for each $x$. Moreover, $L''(v)=4$ and $L''(x)=3$ if $x$ is adjacent to $v$ except to $x_1$  and $x_2$. There is a color $a\in L''(v)-(L''(x_1)\cap L''(x_2))$, WLOG, we let $a\in L''(x_1)$. We color $v$ with a color $a$. Moreover, $|L''(x_1)-\{a\}|=1$ and $|L''(x)-\{a\}|\geq2$ for the others vertices. Thus It is easy that $Z$ has an $L''$-coloring, a contradiction.
\end{pr}
\begin{thm}\label{thmc}  If each vertex of $W_5$ in $G$ is a $5^-$-vertex, then  $W_5$ has at least three $5$-vertices. 
\end{thm}
\begin{pr} Let $V(W_5)=x_1,x_2,x_3,x_4,x_5$ where $x_5$ be a hub and $d(x_i)\leq5$ for each $i$. WLOG, we suppose that $d(x_1)$, $d(x_2)=5$ and the others vertices are $4$-vertices. if $x_1$ is not adjacent to $x_2$, then it 's a contradiction by Lemma \ref{lema}. Let $L$ be a $4$-list assignment of $G$. 
By minimality of $G,$ the graph $G -V(W_5)$ has an $L'$-coloring. If we show that $W_5$ has an $L''$-coloring for any  $|L''(x_i)|=4-|N(x_i)-W_5|)$, then $G$ is $4$-choosable, a contradiction. 
Then, we have  $L''(x_1)=2$, $L''(x_2)=2$, $L''(x_3)=3$, $L''(x_4)=3$, and $L''(x_5)=4$. If there is a color $a$ from $L''(x_5)$ such that $a\notin L''(x_1)\cup L''(x_2)$, then $|L''(x_i)-\{a\}|\geq 2$ for each $i\neq 5.$ It 's complete since a $4$-cycle is $2$-choosable. 
Otherwise, we have $L''(x_1)\cap L''(x_2)=\emptyset$. WLOG, we let $L''(x_1)=\{a,b\}$. First, we choose $a$ from $L''(x_5)$. 
Thus there is an $L''$-coloring $f$ where $f(x_5)=a$, $f(x_1)=b$, $f(x_3)\in  L''(x_3)-\{a,b\}$, 
$f(x_4)\in L''(x_4)-\{a, f(x_3)\}$, and $f(x_2)\in L''(x_2)-\{f(x_4)\}$, a contradiction. Thus $W_5$ has at least three $5$-vertices. 
\end{pr}

\section{Proof of Theorem \ref{main2}} 

\indent Embed a minimal counterexample graph $G$ into the plane. 
Let the initial charge of a vertex $u$ in $G$ be $\mu(u)=2d(u)-6$ 
and the initial charge of a face $f$ in $G$ be $\mu(f)=d(f)-6$. 
Then by Euler's formula $|V(G)|-|E(G)|+|F(G)|=2$ and by the Handshaking lemma, we have
$$\displaystyle\sum_{u\in V(G)}\mu(u)+\displaystyle\sum_{f\in F(G)}\mu(f)=-12.$$ 
\indent Now we design the discharging rule transferring charge 
from one element to another to provide a new charge $\mu^*(x)$ 
for all $x\in V(G)\cup F(G).$ The total of new charges remains $-12$.  
If the final charge  $\mu^*(x)\geq 0$ for all $x\in V(G)\cup F(G)$, 
then we get a contradiction and the proof is completed.\\ 

Before we establish a discharging rule, some definitions are required. 
A cluster of three $3$-faces isomorphic to a graph 
consist of a vertex set of five elements, namely $\{u, v, w, x, y\},$ 
and an edge set $\{xy, xu, xv, yv, yw, uv, vw\}$ 
is  called a \textit{trio}. 
A vertex that is not in any trio is called a \textit{good} vertex. 
We call a vertex $s$ on a face $f$ in a trio 
a \textit{bad} vertex of $f$ if $f$ is the only $3$-cycle 
containing $s$ on that trio, 
a \textit{worst} vertex of $f$ if $s$ is a vertex of  all three $3$-cycles in a trio, 
otherwise $s$ is called a \textit{worse} vertex of $f.$  
We call a face $f$ is a \textit{bad} (\textit{worse},
 or \textit{worst}, respectively) face  
of a vertex $v$ if $v$ is a bad (worse, or worst, respectively) vertex of  $f.$ 
Note that each external vertex of $W_5$ formed by four $3$-faces 
is a worse vertex of some trio.

\indent Let $w(v \rightarrow f)$ be the charge transfered from 
a vertex $v$ to an incident face $f.$ 
The discharging rules are as follows.\\
\textbf{(R1)} Let $f$ be a $3$-face that is not adjacent to the others $3$-faces.\\
\indent \textbf{(R1.1)} For a $4$-vertex $v$,\\
\indent \indent 
$
 w(v \rightarrow f) =
  \begin{cases}
     0.6,        & \text{if } v \text{ is flaw where } f \text{ is a } (4,5^+,5^+)\text{-face},\\
   0.8,        & \text{if } v \text{ is flaw where } f \text{ is a } (4,4,5^+)\text{-face},\\
   1,        & \text{otherwise.}
  \end{cases}
$\\
\indent \textbf{(R1.2)} For a $5^+$-vertex $v$,\\ 
\indent \indent 
$
 w(v \rightarrow f) =
  \begin{cases}
   1.4,        & \text{if } f  \text{ is a } (4,4,5^+)\text{-face where each } 4\text{-vertex is flaw},\\
     1.2,       & \text{if } f  \text{ is a } (4,4^+,5^+)\text{-face where exactly one } 4\text{-vertex is flaw},\\
   1,        & \text{otherwise.}
  \end{cases}
$\\
\textbf{(R2)} Let $f$ be a $3$-face that is adjacent to the others $3$-faces.\\
\indent \textbf{(R2.1)} For a $4$-vertex $v$,\\ 
\indent \indent 
$
 w(v \rightarrow f) =
  \begin{cases}
   0.5,        & \text{if } v \text{ is a hub of } W_5, \\
   1,        & \text{if } f \text{ is a good, bad, or worse face of } v,\\
   2/3,        & \text{if } f \text{ is a worst face of } v.
  \end{cases}
$\\
\indent \textbf{(R2.2)} For a $5$-vertex $v$,\\
\indent \indent 
$
 w(v \rightarrow f) =
  \begin{cases}
   1,        & \text{if } f \text{ is a good or worst face of } v,\\
      1.5,        & \text{if } f \text{ is a bad face of } v,\\
   1.25,        & \text{if } f \text{ is a worse face of } v.
  \end{cases}
$\\
\indent \textbf{(R2.3)} For a $6^+$-vertex $v$,\\
\indent \indent 
$
 w(v \rightarrow f) =
  \begin{cases}
   1,        & \text{if } f \text{ is a good or worst face of } v,\\
   1.5,        & \text{if } f \text{ is a bad or worse face of } v.
  \end{cases}
$\\
\textbf{(R3)} Let $f$ be a $4$-face.\\
\indent \textbf{(R3.1)} For a $4$-vertex $v$, 
\indent \indent 
$w(v \rightarrow f) = 1/3$.\\
\indent \textbf{(R3.2)} For a $5^+$-vertex $v$,\\
\indent \indent 
$
w(v \rightarrow f) =
  \begin{cases}
   1,        & \text{if } f \text{ is a } (4,4,4,5^+)\text{-face},\\
   2/3,        & \text{otherwise. } 
  \end{cases}
$\\

\textbf{(R4)} Let $f$ be a $5$-face.\\
\indent \textbf{(R4.1)} For a $4$-vertex $v$,\\ 
\indent \indent
$
 w(v \rightarrow f) =
  \begin{cases}
     0,        & \text{if } v  \text{ is a flaw and each adjacent vertex of } v \text{ is  a } 4\text{-vertex},\\
     0.1,        & \text{if } v  \text{ is a flaw and one adjacent vertex of } v \text{ is  a } 5^+\text{-vertex},\\
   0.2,        & \text{if } v  \text{ is a flaw and at least two adjacent vertices of } v \text{ are } 5^+\text{-vertices},\\
       1/3,        & \text{if } v \text{ is not a flaw } 4\text{-vertex.}
  \end{cases}
$\\
\indent \textbf{(R4.2)} For a $5$-vertex $v$,\\
\indent \indent 
$
 w(v \rightarrow f) =
  \begin{cases}
   0.4,        & \text{if } f \text{ is a } (4,4,4,5,5)\text{-face where both incident } 5\text{-vertices are adjacent},\\
     1/3,        & \text{if } f \text{ is a } (4,4,4,4,5)\text{-face},\\
       0.3,        & \text{otherwise.}
  \end{cases}
$\\
\indent \textbf{(R4.3)} For a $6$-vertex $v$,\\ 
\indent \indent 
$
 w(v \rightarrow f) =
  \begin{cases}
   0.8,        & \text{if } f \text{ is a } (4,4,4,4,6)\text{-face},\\
       0.4,        & \text{otherwise.}
  \end{cases}
$\\
\indent \textbf{(R4.4)} For a $7^+$-vertex $v$,
$f$ $w(v \rightarrow f) = 0.8$.\\
\textbf{(R5)} Let $f$ be a $7^+$-face and $g$ be a $3$-face.  
If $g$ and other three $3$-faces form $W_5$ and 
$f$ shares exactly one edge with $g,$ let $w(f \rightarrow g) = 1/8.$\\
\textbf{(R6)} After (R1) to (R6),  
redistribute the total of charges of $3$-faces in 
the same cluster of adjacent $3$-faces (trio or $W_5$) equally among 
its $3$-faces.\\

\indent It remains to show that resulting $\mu^*(x)\geq 0$ 
for all $x\in V(G)\cup F(G)$.\\
\indent\textit{CASE 1: } Consider  a $4$-vertex $v.$\\
\indent We use  (R1.1), (R2.1), (R3.1), and (4.1) to prove this case.\\ 
\textit{SUBCASE 1.1: } Let $v$ be flaw $4$-vertex.\\ \indent Then $v$ is a $(3,3,5,5^+)$-vertex. If each adjacent vertex of $v$ is a $4$-vertex, then  we obtain $\mu^*(v) = \mu(v) - (2\cdot 1)= 0$. If $v$ is adjacent to exactly one $5^+$-vertex, then we obtain $\mu^*(v) = \mu(v) - (1+0.8+2\cdot (0.1)) \geq 0$. If $v$ is adjacent to at least two $5^+$-vertices, the we obtain $\mu^*(v) = \mu(v) - (2\cdot(0.8)+2\cdot (0.2)) \geq 0$ for two $5^+$-vertices are not in the same $3$-face. Otherwise, we obtain $\mu^*(v) = \mu(v) - (1+0.6+2\cdot (0.2)) \geq 0$.\\
\textit{SUBCASE 1.2: } Let $v$ be not  flaw $4$-vertex.\\
\indent  If $v$ is incident at most one $3$-face, then we obtain $\mu^*(v)\geq \mu(v) - (1 + 3\cdot (1/3)) \geq 0$. If $v$ is  incident to two $3$-faces, then $v$ is a $(3,3,6^+,6^+)$-vertex. Thus we obtain $\mu^*(v) = \mu(v) - (2\cdot 1)= 0$.  If $v$ is incident to  three $3$-faces, 
then $v$ is a worst vertex of these faces 
and its remaining incident face is a $6^+$-face.  
Thus we obtain  $\mu^*(v) = \mu(v) - (3\cdot (2/3)) = 0.$ If $v$ is incident to  four $3$-faces, 
then $v$ is a hub of $W_5$. 
Thus $\mu^*(v)= \mu(v) - 4\cdot (1/2) = 0.$\\
\indent \textit{CASE 2: } Consider a $5$-vertex $v.$\\ 
\textit{SUBCASE  2.1:} Let $v$ be incident to an adjacent triangle or a bad face. \\
\indent Then  $v$ is incident to at least two $6^+$-faces by Proposition \ref{main1}(1). Additionally, a $5$-vertex $v$ has at most two bad faces. We use  (R2.2), (R3.2), and (R4.2) to prove the following cases.\\
\indent  We obtain $\mu^*(v) = \mu(v) - (3\cdot(1.25)) >0$ if there is no any bad face, $\mu^*(v) = \mu(v) - (1.5+2\cdot(1.25)) \geq 0$ if there is one bad face, and $\mu^*(v) = \mu(v) - 2\cdot(1.5) >0$ if there are two bad faces.\\
\textit{SUBCASE  2.2:} A vertex $v$ has neither adjacent triangles nor  bad faces.\\
\indent Then $v$ is incident at most two $3$-faces. We use  (R1.2), (R3.2), and (R4.2) to prove the following cases. \\
\indent Let $v$ be incident to at least one $6^+$-face.   
If $v$ is not incident to any $3$-face, then $\mu^*(v) = \mu(v) - (4\cdot(1)) \geq 0$. 
If $v$ is incident exactly one $3$-face, then $\mu^*(v) = \mu(v) - (1.4+2\cdot(1)+0.4) >0$. 
If $v$ is incident to two $3$-faces, then $v$ is incident to at most one $4$-face by Proposition \ref{main1}(2). Then we obtain $\mu^*(v) = \mu(v) - (1.4+1.2+1+0.4) >0$.\\ 
\indent Next, a $5$-vertex $v$ is not incident to any $6^+$-face.\\
\indent Let $v$ be not incident to any $3$-face. 
If  $v$ is  a $(4^+,4^+,4^+,5,5)$-vertex, then $\mu^*(v) = \mu(v) - (3\cdot(1)+2\cdot(0.4))>0$. 
If  $v$ is  a $(4,4,4,4,5)$-vertex, then at least two incident $4$-faces of $v$ are $(4^+,4^+,5,5^+)$-faces by Theorem \ref{thma}. 
Thus we obtain 
 $\mu^*(v) = \mu(v) - (2\cdot(1)+2\cdot(2/3))+0.4>0$.
If  $v$ is a $(4,4,4,4,4)$-vertex, then then at least three incident $4$-faces of $v$ are $(4^+,4^+,5,5^+)$-faces by Theorem \ref{thma}. 
Thus we obtain 
 $\mu^*(v) = \mu(v) - (2\cdot(1)+3\cdot(2/3))>0$.\\ 
\indent Let $v$ be  incident to one $3$-face. 
If  $v$ is a $(3,4^+,5,5,5)$-vertex, then $\mu^*(v) = \mu(v) - (1.4+1+3\cdot(0.4))\geq 0.4$. 
Let $v$ be a $(3,4,4,5,5)$-vertex with two $4$-faces, $f_1$ and $f_2$. If $f_1$ is adjacent to  $f_2$, then either $f_1$ or $f_2$ is a $(4^+,4^+,5,5^+)$-face by Lemma \ref{lema}. 
We obtain $\mu^*(v) = \mu(v) - (1.4+1+2/3+2\cdot(0.4))>0$. 
If $f_1$ is not adjacent to  $f_2$, then an incident $3$-face of $v$ is adjacent to $f_1$ and $f_2$ by Proposition \ref{main1}(3). Then we obtain 
$\mu^*(v) = \mu(v) - (3\cdot(1)+2\cdot(0.4))\geq 0.2.$\\
\indent  
Let $v$ be incident to two $3$-faces. 
By Proposition \ref{main1}(2), $v$ is not incident to at least two $4$-faces and by  Proposition \ref{main1}(3), $v$ is not incident to one $4$-face. 
Thus $v$ is  a $(3,3,5,5,5)$-vertex and we obtain $\mu^*(v) = \mu(v) - (2\cdot(1.4)+3\cdot(0.4))\geq0.$ \\
\indent \textit{CASE 3: } Consider a $6$-vertex $v.$\\ 
\textit{SUBCASE  3.1:} Let $v$ be incident to an adjacent triangle or a bad face. \\
\indent Then  $v$ is incident to at least two $6^+$-faces by Proposition \ref{main1}(1). Thus we obtain $\mu^*(v) = \mu(v) - 4\cdot(1.5)\geq0$ since $v$ sends charge at most $1.5$ to each incident face by (R1.2), (R2.3), (R3.2), and (R4.3).\\
\textit{SUBCASE  3.2:} A vertex $v$ has neither adjacent triangles nor  bad faces.\\ 
\indent Then $v$ is incident at most three $3$-faces. We use  (R1.2), (R3.2), and (R4.3) to prove the following cases.\\
\indent Let $v$ be incident to at least one $6^+$-face. Then 
 $\mu^*(v) = \mu(v) - (2\cdot(1.4)+3\cdot(1)) >0$ if there are at most two $3$-faces and we obtain  $\mu^*(v) = \mu(v) - (3\cdot(1.4)+2\cdot(0.8)) >0$ if there are three $3$-faces.\\  
 \indent Next, a $6$-vertex $v$ is not incident to any $6^+$-face.\\ 
 \indent Let $v$ be not incident to any $3$-face. 
 Then $\mu^*(v) = \mu(v) - (6\cdot(1)) \geq 0.$\\ 
\indent Let $v$ be incident to one $3$-face. 
By Proposition \ref{main1}(2), $v$ is incident to at most three $4$-faces. 
If $v$ is a $(3,4,4,4,5,5)$-vertex, then $\mu^*(v) = \mu(v) - (1.4+3\cdot(1)+2\cdot(0.8))\geq0$.\\  
\indent Let $v$ be incident to two $3$-faces.
By Proposition \ref{main1}(2) and \ref{main1}(3), $v$ is incident to at most one $4$-face.   
If  $v$ is a $(3,3,4,5,5,5)$-vertex, then an incident  $4$-face of $v$ and an incident $3$-faces of $v$ are not adjacent by Proposition \ref{main1}(3). 
By Theorem \ref{thmb},  at least two incident faces of $v$ are incident at least two $5^+$-vertices. Then one of incident $5$-face of $v$ is incident to two $5^+$-vertices. We obtain  $\mu^*(v) = \mu(v) - (2\cdot(1.4)+1+2\cdot(0.8)+0.4)>0.$ 
If $v$ is a $(3,3,5,5,5,5)$-vertex, then $\mu^*(v) = \mu(v) - (2\cdot(1.4)+4\cdot(0.8))\geq0.$\\ 
\indent Let $v$ be incident to three $3$-faces.
By Proposition \ref{main1}(2), $v$ is not  incident to any $4$-face. Thus 
 $v$ is a $(3,3,3,5,5,5)$-vertex. 
By Theorem \ref{thmb},  at least two of incident faces $f_1$ and $f_2$ of $v$ are incident  to 
at least two $5^+$-vertices. Then we obtain $\mu^*(v) = \mu(v) - (3\cdot(1.4)+0.8+2\cdot(0.4))\geq0$ if $f_1$ and $f_2$ are $5^+$-faces  and $\mu^*(v) = \mu(v) - (2\cdot(1.4)+1.2+2\cdot(0.8)+0.4)\geq0$ if  a $3$-face is either $f_1$ or $f_2$. \\
\indent \textit{CASE 4: } Consider a $k$-vertex $v$ with $k\geq 7.$\\
\textit{SUBCASE 4.1: }Let $v$ be incident to an adjacent triangle or a bad face. \\
\indent Then  $v$ is incident to at least two $6^+$-faces by Proposition \ref{main1}(1). Then $v$ sends charge to  at most $k-2$ faces. We obtain $\mu^*(v) = \mu(v) - (k-2)\cdot(1.5)>0$ since $(2\cdot(k)-6)/k-2\geq1.5$ for $k\geq6$ by (1.2), (R2.3), (R3.2), and (R4.4).\\
\textit{SUBCASE:  4.2} A vertex $v$ has neither adjacent triangles nor  bad faces.\\
\indent We use  (R1.2), (R3.2), and (R4.4) to prove the following cases.\\
\indent Let $v$ be a $7$-vertex. 
If $v$ is incident to at most two $3$-faces, then  
$\mu^*(v) = \mu(v) - (2\cdot(1.4)+5\cdot(1))\geq0.2$.  
If $v$ is incident to three $3$-faces, then $v$ is incident to at most one $4$-face by Proposition \ref{main1}(2). Thus we obtain 
$\mu^*(v) = \mu(v) - (3\cdot(1.4)+1\cdot(1)+3\cdot(0.8))\geq0.4.$\\ 
\indent Let $v$ be a $8^+$-vertex. Then $v$ is incident to at most $d(v)/2$ $3$-faces. $\mu^*(v) \geq0.4$ since $1.4(d(v)/2)+d(v)/2\leq2d(v)-6$ for $d(v)\geq8$. 
Thus $\mu^*(v) \geq0.4$ for each vertex $v.$\\
\indent It is clear  that $\mu^*(f)=\mu(f) \geq 0$ for each $f$ is a $6^+$-face  since 
$(1/8)d(f)<d(f)-6$ for $d(f)\geq7.$\\
\indent \textit{CASE 5: }  Consider a $3$-face $f$ that is not adjacent to the others $3$-faces.\\ 
\indent We use 
(R1.1) and (R1.2) to prove the following cases.\\   
\indent Let $f$ be not incident to each flaw $4$-vertex. Then $\mu^*(f) = \mu(f) +(3\cdot(1))=0$.\\ 
 \indent Next at least one incident $4$-vertex of $f$ is flaw.  
If $f$ is a $(4,4,5^+)$-face, then $\mu^*(f) = \mu(f) +(1.4+2\cdot(0.8))=0$ when both incident $4$-vertices of $f$ are flaw and $\mu^*(v) = \mu(v) +(1.2+1+0.8)=0$ when exactly one of incident $4$-vertex of $v$ is flaw. 
If $f$ is a $(4,5^+,5^+)$-face, then $\mu^*(f) = \mu(f) +(2\cdot(1.2)+0.6)=0$.\\
\indent \textit{CASE 6: } Consider a $3$-face $f$ that is  adjacent to the others $3$-faces.\\
\indent We use 
(R2.1), (R2.2), (R2.3), (R5), and (R6) to prove the following cases.\\
\textit{SUBCASE 6.1:} If $f$ is not in a trio, then $\mu^*(f) = \mu(f) +(3\cdot(1))=0$ since each incident vertex sends charge at least one to $f$.\\
\textit{SUBCASE 6.2:} Let $f$ be in a trio.\\ 
\indent Let $f_1,f_2,$ and $f_3$ be $3$-faces in a same trio $T.$ 
Define $\mu(T) := \mu(f_1)+\mu (f_2)+ \mu(f_3) = -9$ and 
 $\mu^*(T) := \mu^*(f_1)+\mu^* (f_2)+ \mu^*(f_3)$ by (R6).\\
 \indent If $f$ is in a trio that a worst vertex is not a $4$-vertex, then each $3$-face of trio $h$ that  $\mu^*(h) = \mu(h) +(3\cdot(1)) \geq 0$.\\
\indent If $T$ is a trio that a worst vertex is a $4$-vertex, then there are many cases as follows.\\ 
\indent If each worse vertex is a $4$-vertex, then two bad vertices are  $5^+$-vertices by Lemma \ref{lema}. 
Then  $\mu^*(T) = -9 +3\cdot(2/3)+2\cdot(1.5)+4\cdot(1) \geq 0$.\\
\indent If one of worse vertex is a $5$-vertex, then  either the other worse vertex or at least one bad vertex is a $5^+$-vertex by Lemma \ref{lema}. 
Then $\mu^*(T) = -9 +3\cdot(2/3)+4\cdot(1.25)+2\cdot(1)  \geq 0$ or 
$\mu^*(T) = -9 +3\cdot(2/3)+2\cdot(1.25)+1.5+3\cdot(1) \geq 0$, respectively.\\ 
\indent If a worse vertex is a $6^+$-vertex, then $\mu^*(T) = -9 +3\cdot(2/3)+2\cdot(1.5)+4\cdot(1) \geq 0$. 
\\
\textit{SUBCASE 6.3:} Let $f$ be in $W_5$.\\ 
\indent If each vertex of $W_5$ is not a $6^+$-vertex, then at least three vertices are $5$-vertices by Theorem \ref{thmc}. Thus  we obtain $\mu^*(W_5) = -12  +6\cdot(1.25)+2\cdot(1)+4\cdot(0.5)+4\cdot(1/8)\geq 0.$\\
\indent  If  exactly one vertex of $W_5$ is a $6^+$-vertex, then one of the others vertices is a $5^+$-vertex by Lemma \ref{lema}. Thus we  obtain $\mu^*(W_5) = -12  +2\cdot(1.5)+2\cdot(1.25)+4\cdot(1)+4\cdot(0.5)+4\cdot(1/8)\geq 0.$\\
\indent  If  at least two vertices of $W_5$ are  $6^+$-vertices, then  we obtain $\mu^*(W_5) = -12  +4\cdot(1.5)+4\cdot(1)+4\cdot(0.5)+4\cdot(1/8)> 0$. 
\\
\indent \textit{CASE 7: } Consider a $4$-face $f.$\\ 
\indent Then every vertex on $f$ has degree at least $4$ and one 
of them has degree at least $5.$ 
If $f$ is a $(4,4,4,5^+)$-face, 
then $\mu^*(f)\geq \mu(f) + 3(1/3)+1 =0$ and      
we obtain  $\mu^*(f) \geq \mu(f) +2(1/3)+2 \cdot (2/3) \geq 0$ if  $f$ is a $(4^+,4^+,5^+,5^+)$-face. \\
\indent \textit{CASE 8: } Consider a $5$-face $f.$\\
 \indent We use 
(R4.1), (R4.2), (R4.3), and (R4.4)  to prove the following cases.\\ 
\textit{SUBCASE 8.1: } Let $f$ be incident to at least three $5^+$-vertices.\\
 \indent Then, 
each incident  $4$-vertex of $f$ is adjacent to at least one $5^+$-vertex. Thus we obtain  $\mu^*(f)=\mu(f) +3\cdot(0.3)+2\cdot(0.1)  \geq 0$.\\  
\textit{SUBCASE 8.2: } Let $f$ be incident to two $5^+$-vertices $x$ and $y$.\\
\indent 
If  $x$ is not adjacent to $y$, then each incident  $4$-vertex of $f$ is adjacent to at least one $5^+$-vertex. Moreover, one of these is adjacent to at least two $5^+$-vertices. Thus we obtain   $\mu^*(f)=\mu(f)+ 2\cdot(0.3)+2\cdot(0.1)+0.2 \geq 0.$\\  
\indent If  $x$ is adjacent to $y$, then two of incident $4$-vertices of $f$ is adjacent to at least one $5^+$-vertex. Thus we obtain  $\mu^*(f)=\mu(f)+ 2\cdot(0.4)+2\cdot(0.1)  \geq 0.$.\\ 
\textit{SUBCASE 8.3: }Let $f$ be incident to at most one $5^+$-vertex.\\ 
\indent If $f$ is a $(4,4,4,4,6^+)$-face, then two of incident $4$-vertices of $f$ are adjacent to at least one $6^+$-vertex. Thus we obtain  $\mu^*(f)=\mu(f)+ 0.8+2\cdot(0.1)  \geq 0.$\\ 
\indent Now, it remains to show $f$ is a $(4,4,4,4,5^-)$-face.\\ 
\indent If $f$ is adjacent to at most one $4^+$-face,  then 
$f$ is a poor $5$-face. Thus each incident flaw $4$-vertex of $f$ is adjacent to at least two $5^+$-vertices  by Lemma \ref{poor}. Thus we obtain 
$\mu^*(f)=\mu(f) +5\cdot(0.2)  \geq 0.$\\ 
\indent If $f$ is adjacent to at least two $4^+$-faces, then at least three incident  vertices of $f$ are not flaw $4$-vertices. Thus we obtain $\mu^*(f)=\mu(f) +3\cdot(1/3)  \geq 0.$\\
\indent This completes the proof. 

\section{Acknowledgments}
\indent The first author is supported by Development and Promotion 
of Science and Technology Talents Project (DPST).


\newpage

\begin{thebibliography}{99}

\bibitem{app1} K. Appel, W. Haken, 
Every planar map is four colorable. I. 
Discharging, Illinois J. Math. 21(3)(1977) 429-490.

\bibitem{app2} K. Appel, W. Haken, J. Koch, 
Every planar map is four colorable. II. 
Reducibility, Illinois J. Math. 21(3)(1977) 491-561.

\bibitem{Bo}
O.V. Borodin, A.O. Ivanova, 
Planar graphs without triangular $4$-cycles are $4$-choosable, 
Sib. \`{E}lektron. Mat. Rep. 5(2008) 75-79.

\bibitem{Cheng}
P.P. Cheng, M. Chen, Y.Q. Wang, 
Planar graphs without $4$-cycles adjacent to triangles are $4$-choosable, 
Discrete Math. 339(2016) 3052-3057.

\bibitem{Erdos}
P. Erd\H os, A.L. Rubin, H. Taylor, 
Choosability in graphs, in: 
Proceedings, West Coast Conference on Combinatorics, 
Graph Theory and Computing, Arcata,
CA., Sept. 5-7, in: Congr. Numer., vol. 26, 1979.

\bibitem{Far}
B. Farzad, Planar graphs without $7$-cycles are $4$-choosable, 
SIAM J. Discrete Math. 23(2009) 1179-1199.

\bibitem{Fi}
G. Fijav\v z, M. Juvan, B. Mohar, R. \v Skrekovski, 
Planar graphs without cycles of specific lengths, European J. Combin. 23(2002) 377-388.

\bibitem{Gut}
S. Gutner, The complexity of planar graph choosability, 
Discrete Math. 159(1996) 119-130.

\bibitem{Hu}
D.Q. Hu, J.L. Wu, Planar graphs without intersecting $5$-cycles are $4$-choosable, 
Discrete Math. 340(2017) 1788-1792. 

%\bibitem{Lam1}
%P.C.B. Lam, W.C. Shiu, B.G. Xu, On structure of some plane graphs with applications to choosability, J. Combin. Theory Ser. B 82(2001) 285-296.

\bibitem{Lam2} 
P.C.B. Lam, B. Xu, J. Liu, The $4$-choosability of plane graphs without $4$-cycles, J. Combin. Theory Ser. B 76(1999) 117-126.

\bibitem{Tho}
C. Thomassen, Every planar graph is $5$-choosable, 
J. Combin. Theory Ser. B 62(1994) 180-181.

\bibitem{Vizing} 
V.G. Vizing, Vertex colorings with given colors, 
Metody Diskret. Analiz. 29(1976) 3-10 (in Russian).

\bibitem{Vo1}
M. Voigt, List colourings of planar graphs, Discrete Math. 120(1993) 215-219.

\bibitem{Xu} 
R. Xu, J.L. Wu, 
A sufficient condition for a planar graph to be $4$-choosable, 
Discrete App. Math. 224(2017)120-122.

\bibitem{Wang1}
W. Wang, K.W. Lih, 
Choosability and edge choosability of planar graphs without five cycles, 
Appl. Math. Lett. 15(2002) 561-565.

\bibitem{Wang2} 
W. Wang, K.W. Lih, 
Choosability and edge choosability of planar graphs without intersecting triangles, 
SIAM J. Discrete Math. 15(2002) 538-545.



\end{thebibliography}
\end{document}